\documentclass[11pt]{amsart}

\usepackage{amscd}
\usepackage{amssymb}
\usepackage[latin1]{inputenc}

\def\Z{\mathbb{Z}}

\def\ZZ{\Z \oplus \Z}

\def\g{\gamma}

\def\l{\lambda}

\def\s{\sigma}
\def\w{\omega}
\def\G{\Gamma}

\def\cal{\mathcal}

\def\ol{\overline}

\newtheorem{thm}{Theorem}[section]

\newtheorem{lem}{Lemma}[section]

\title{Conjugacy problem in groups of non-oriented geometrizable
3-manifolds}

\begin{document}

\maketitle
\begin{center}
{\sc Jean-Philippe PR\' EAUX}\footnote[1]{Centre de recherche de
l'Ecole de l'air, Ecole de l'air, F-13661 Salon de Provence air}\
\footnote[2]{Centre de Math\'ematiques et d'informatique,
Universit\'e de Provence, 39 rue F.Joliot-Curie, F-13453 marseille
cedex 13\\
\indent {\it E-mail :} \ preaux@cmi.univ-mrs.fr\\
\indent {\sl Keywords :} \ Conjugacy problem, fundamental group,
non-oriented 3-manifold, Thurston's geometrization conjecture.\\
\indent{\sl AMS subject classification :} \ 57M05, 20F10.}
\end{center}

\begin{abstract}
We have proved in \cite{cp3mg} that fundamental groups of oriented
geometrizable  3-manifolds have a solvable conjugacy problem. We
now focus on groups of non-oriented geometrizable 3-manifolds in
order to conclude that all groups of  geometrizable 3-manifolds
have a solvable conjugacy problem.
\end{abstract}

\section*{Introduction}
Since enonced by M.Dehn in the early 1910's the word problem and
conjugacy  problem in finitely presented groups have become
fundamental problems in combinatorial group theory. Following the
work of Novikov \cite{novikov} and further authors on their
general unsolvability, it has become fairly natural to ask for any
finitely presented group whether it admits a solution or not. For
example in \cite{dehn1,dehn2,dehn3}, Dehn has solved those
 problems for fundamental groups of surface ; its motivation was
topological.

 Given a finite
presentation of a group $G$, a solution to the
 {\sl word problem}  is an algorithm which given
  $\w,\w'\in G$ (as a couple of words on the generators),
   decides whether $\w=\w'$ in $G$ or not. A solution to the {\sl conjugacy problem}
    is an algorithm which given $u,v\in G$ decides whether
     $\exists\ h\in G$ such that $u=hvh^{-1}$ in $G$ or not.

It turns out that existence of solutions does not depend of the
finite presentation involved. Hence existence of a solution in $G$
to any of these problems  only depends on the isomorphism class of
$G$. We say that $G$ has a {\sl solvable} word problem (resp.
conjugacy problem) if  $G$ admits a solution to the word problem
(resp. conjugacy problem).

By a 3-manifold we mean a connected compact manifold with boundary
of dimension 3 ; a 3-manifold may be oriented or not. We work in
the $PL$ category  ; according to the hauptvermutung and Moise
theorem this is not restrictive. Following the work of Thurston
(cf. \cite{th}) an oriented 3-manifold $M$ is {\sl geometrizable}
if the pieces obtained in its canonical topological decomposition
 (roughly speaking along essential spheres, discs and tori) have interiors
which admit complete locally homogeneous riemanian metrics. A non
oriented 3-manifold is said to be geometrizable, if its
orientation cover  is geometrizable.

In the class of fundamental groups of geometrizable 3-manifolds,
the word problem is  known to be solvable, following early work of
Waldhausen (\cite{wald}) as well as more recent work of Epstein
and Thurston  on automatic group theory (cf. \cite{epstein}). We
make use of this result in our proof.

We have considered in \cite{cp3mg} the conjugacy problem in groups
of oriented geometrizable 3-manifolds. We are now  concerned with
the non-oriented case to provide a conclusion concerning the
conjugacy problem in the whole class of groups of geometrizable
3-manifolds.

\section{Statement of the result}
 This work is entirely devoted to the proof of the following result :
\begin{thm}\label{main}
Fundamental groups of non-oriented geometrizable 3-manifolds have solvable conjugacy problem.
\end{thm}
It  does not follow as a  consequence of existence of a solution
in the oriented case (cf. \cite{cp3mg}), since D.Collins and
C.Miller have shown that the conjugacy problem can be unsolvable
in a group even when solvable in an index 2 subgroup
(\cite{collins}).
 Nevertheless
our strategy will consist essentially in reducing to the conjugacy
problem in the oriented case. Together with a solution in the
oriented case (cf. \cite{cp3mg}), one obtains
\begin{thm}
Fundamental groups of geometrizable 3-manifolds have solvable conjugacy problem.
\end{thm}
\noindent as well as the following corollaries. The former one is
a topological rephrasing : a loop $\g$ in $M$ defines a conjugacy
class in $\pi_1(M)$ which only depends on its free homotopy class.

\begin{thm}
Given a geometrizable 3-manifold $M$ there exists an algorithm
which decides for any couple of loops $\g,\, \g'$ in $M$  whether
$\g$ and $\g'$ are freely homotopic.
\end{thm}

Concerning decision problems relative to boundary subgroups, one has :

\begin{thm}
\label{t2} Let $M$ a  geometrizable 3-manifold $M$, $F$ any
compact connected surface lying in $\partial M$ and $G=\pi_1(M)$,
$H=i_*(\pi_1(F))$ ; there exists an algorithm which  decides for
any $g\in G$ whether $g\in H$ (resp. $g$ is conjugate to an
element of $H$).
\end{thm}

And one has the topological rephrasing :

\begin{thm}
 Given any 3-manifold $M$ with $\partial M\not= 0$, and $F$ any connected surface in $\partial M$,
  one can decide for any loop $\g$ (resp. $*$-based loop, $*\in F$) in $M$ whether
  up to
   homotopy (resp. homotopy with $*$ fixed) $\g$ lies in $F$.
\end{thm}

\section{Some decision results in extensions by $\Z_2$}

Let $G$ a group and $u,v,h\in G$. Once $u=hvh^{-1}$ we shall use
the notation $u=v^h$ or $u\sim v$. We denote $Z_G(v)=\{u\in G\ |\
uv=vu\}$ the centralizer of $v$ in $G$ and $C_G(u,v)=\{h\in G\ |\
u=v^h\}$. The subset $C_G(u,v)$ of $G$ is either empty (when
$u\not\sim v$) or $h.Z_G(v)$ for some $h\in G$ such that $u=v^h$.

\begin{lem}
\label{ld1}
 Let $G$ be a group and $H$  an index 2 subgroup of $G$ with solvable conjugacy problem.
 Given any couple of elements $u,v\in H$ one can decide whether $u$ and $v$ are conjugate in $G$.
\end{lem}
\noindent {\sl Proof.} Given a set of representative $a_0=1,a_1$
of $H /G$, in order to decide whether $u,v\in H$ are conjugate in
$G$ it suffices to decide whether $u$ is conjugate in $H$ to any
of the $a_iva_i^{-1}$ for $i=0,1$.\hfill $\square$

\begin{lem}
\label{ld2} Let $G$ be a group and $H$ be an index 2 subgroup of
$G$ ; suppose that $H$ and $G$ have respectively solvable
conjugacy problem and solvable word problem. Let $v\in G\setminus
H$ such that
 $Z_G(v)=Z_G(v^2)$ ; then  on can decide for any $u\in G$ whether $u$ and $v$ are conjugate in $G$.
\end{lem}

\noindent {\sl Proof.} Let $v\in G$ be as above, and suppose one
wants to decide for some given $u\in G$ whether $u \sim v$ in $G$.
  With the preceding lemma one can decide whether $u^2$ and $v^2$ are conjugate in $G$. If not
   then $u$ and $v$ are definitely not conjugate in $G$. So suppose $u^2$ and $v^2$ are conjugate in $G$,
    say $u=k.v.k^{-1}$ for some $k\in G$ that one can find using a solution to the word problem in $G$,
    so that $C_G(u^2,v^2)=k.Z_G(v^2)$.
    Obviously $C_G(u,v)\subset C_G(u^2,v^2)$ and moreover since
    $Z_G(v^2)=Z_G(v)$,
if $C_G(u,v)$ is non empty it must equal $C_G(u^2,v^2)$. Hence to
decide whether $u\sim v$ in $G$ it suffices to decide with the
word problem in $G$ whether $u=v^k$ or not.\hfill$\square$\medskip

We  also remark the following observation.

\begin{lem}
\label{ld3} Let $G$ be a f.g. group with solvable word problem,
and $H$ be a finite index subgroup. Then the generalised word
problem of $H$ in $G$ is solvable.
\end{lem}
\noindent {\sl Proof.} Given a set of representative
$a_0,a_1,\ldots , a_n$ of $G/H$ and a finite set of generators of
$H$, as well as $g\in G$, one can enumerate all elements $h\in H$
and check in $G$ whether $g=a_i.h$ for some $i=1,2,\ldots , n$.
This process must terminate to provide the class of $g$ in
$G/H$.\hfill $\square$

This last algorithm  is far from being efficient, but anyway note
that it will be redondant here ; we only give it as an alternative
algorithm for algebra enclined readers while in our proof it will
replaced by a topological argument.

\section{Checking for centralizers in groups of  oriented  3-manifolds}

Let $M$ be an oriented closed irreducible 3-manifold, and $T$ the
set of essential tori in the JSJ decomposition of $M$ ($T$ may be
empty). When $T$ is non empty, $T$ is well defined up to isotopy,
and the open manifold $M\setminus T$ has connected components
which are open manifolds either homeomorphic to Seifert fibered
spaces or to finite volume hyperbolic manifolds. There exists a
canonical map $r:M\setminus T\longrightarrow M$ which is an
embedding. The {\sl Seifert characteristic submanifold} (cf.
\cite{js}) of $M$ is then
defined to be the compact submanifold of $M$ with connected components constructed as follows :\\
$(i)$ For each Seifert fibered space  component $S_i$ of
$M\setminus T$ cut the ends of $r(S_i)$ and then consider its
adherence in $M$. Such components are called {\sl non trivial
Seifert submanifolds}
 $(ii)$ For each torus $T_i$ in $T$ which is not parallel in $M$ to a boundary component of a non trivial Seifert
 submanifold of $M$
  consider a regular
 neighborhood (homeomorphic to $S^1\times S^1\times I$) of $T$ so that they are all pairwise disjoint and  do not intersect
 non trivial Seifert submanifolds. Such coponents are
 called {\sl trivial Seifert submanifolds}.

  If $T$ is empty the
 Seifert characteristic submanifold is defined to be the empty set if $M$
 is not a Seifert fibered space, and $M$ otherwise.
The Seifert characteristic submanifold is well defined up to
isotopy.

 The embedding of components of the Seifert characteristic submanifold induces, up to conjugacy,
   embeddings of their fundamental
  groups in $\pi_1(M)$. Once embeddings are given, their images
  are called {\sl Seifert subgroups} of $\pi_1(M)$ ; they can be
  isomorphic to $\ZZ$ (when induced from a trivial Seifert
  submanifold) and  called {\sl $\ZZ$ Seifert subgroups} ; otherwise they are fundamental groups of
  Seifert spaces and called {\sl non $\ZZ$ Seifert subgroups}.\\

  Now let $M$ be an oriented 3-manifold whose boundary only consists of
  spheres, and such  that the closed 3-manifold $\hat{M}$ obtained by
  attaching balls to $\partial M$ is irreducible. Depending on a given embedding
  of $M$ into $\hat{M}$, the Seifert
  characteristic submanifold of $M$ is defined to be the intersection of the
  Seifert characteristic submanifold of $\hat{M}$ with $M$.   Its components are
  punctured Seifert submanifolds and are defined up to the
  embedding of $M$. We make use of  the terminology of trivial (resp. non trivial)
   Seifert submanifolds of $M$ for components which come from trivial (resp. non trivial)
   Seifert submanifolds of $\hat{M}$. We also define, up to conjugacy, Seifert
  subgroups of $\pi_1(M)=\pi_1(\hat{M})$.\\

\begin{lem}\label{central}
Let $M$ be an oriented 3-manifold such that $\hat{M}$ is closed
irreducible. Suppose a triangulation of $M$ as well as an
embedding of $M$ into $\hat{M}$ are given.  Then one can
algorithmically construct the Seifert characteristic submanifold
of $M$ and recognize trivial Seifert submanifolds from other
Seifert submanifolds. Given in $\pi_1(M)$ an element $u$ and  a
Seifert subgroup $K$, one can decide whether $u$ is conjugate to
an element of $K$ and, if so, produce a conjugating element.
\end{lem}
\noindent {\sl Proof.} Glue $PL$-balls to $\partial M$ in order to
obtain $\hat{M}$ together with  a triangulation. Then apply the
algorithm of \cite{jt} in order to obtain the JSJ decomposition
$T$ of $\hat{M}$ as well as the Seifert characteristic
submanifold. Note that the algorithm also produces Seifert
invariants of the Seifert submanifolds so that one can recognize
trivial one from the others. Since the embedding of $M$ into
$\hat{M}$ is given it provides the Seifert characteristic
submanifold of $M$. Note that $\hat{M}$ is a $S^1\times
S^1$-bundle over $S^1$ modelled on geometry $Sol$ if and only if
$M\setminus T$ is homeomorphic to the interior of $S^1\times
S^1\times I$, which can be easily checked.

The JSJ decomposition of $\hat{M}$ splits
$\pi_1(M)=\pi_1(\hat{M})$ into a fundamental group of a graph of
group. One can produce such a graph of group as well as a finite
presentation of $\pi_1(M)$ and  embeddings of edge and vertex
groups in $\pi_1(M)$, so that $\ZZ$-Seifert subgroups (resp. non
$\ZZ$ Seifert subgroups) are edge subgroups (resp. vertex
subgroups) of $\pi_1(M)$. In the following we exlusively make use
of results established in \cite{cp3mg}. Given an element
$u\in\pi_1(M)$ one can use corollary 4.1 of \cite{cp3mg} to change
$u$ into a cyclically reduced conjugate according to this
splitting (\S 3.1, \cite{cp3mg}) ; it also provides the
conjugating element. If $\hat{M}$ is a $S^1\times S^1$-bundle over
$S^1$ modelled on geometry $Sol$, then obviously $u$ is conjugate
to $K=\ZZ$ if and only if $u$ lies in the vertex subgroup $K$.
Otherwise, it follows from theorem 3.1, lemma 2.2 and proposition
4.1 of \cite{cp3mg} that $u$ lies in a given vertex subgroup
$G_V$, for some vertex $V$, and that $u$ is conjugate in
$\pi_1(M)$ to an element of $K$ exactly if one of the following
cases occur : $(i)$ $G_V=K$ ; $(ii)$ $K$ is an
    edge subgroup of $G_V$, and $u$ is conjugate in $G_V$ to an element of $K$ ;
    $(iii)$ $K$ is a vertex subgroup $G_{V'}$, for some vertex $V'$ and $u$
    is conjugate in $G_V$
    to an element lying in an edge subgroup $G_e$ for
     some edge $e$ going from $V$ to $V'$. In case $(i)$ or $Sol$ geometry a solution is obvious,
     while cases $(ii)$ and $(iii)$ reduce to
     deciding whether $u$ is conjugate in $G_V$ to a given edge subgroup. One  applies the algorithms of \cite{cp3mg}
      given      in proposition 5.1 (if $G_V$ comes from a Seifert component) or  theorem 6.3 (otherwise)
     to decide so. Note that going into the line of the proofs on sees that the algorithms produce,
     if any, a conjugating element ;
     nevertheless one can also  obviously  find such an element by a naive algorithm.\hfill $\square$ \\

     \noindent{\bf Remark : }
Jaco and Shalen have proved (theorem VI.I.6, \cite{js}) that
whenever $G$ is the group of an oriented closed Haken manifold an
element  $g\in G$ as either a cyclic centralizer or up to
conjugacy its centralizer lies in a Seifert subgroup. If $M$ is an
oriented irreducible non Haken manifold, $M$ must be closed and
homeomorphic to either a Seifert space or an hyperbolic manifold,
and it follows that the centralizers have the same structure as in
the Haken case. If $\partial M$ contains spheres and $\hat{M}$ is
irreducible the same conclusion applies.
 Hence in such cases only
the following can occur :\\
(i) $Z_G(g)=\Z$,\\
(ii) $g$ is conjugate to a $\ZZ$ Seifert subgroup, $Z_G(g)=\ZZ$,\\
(iii) $g$ is conjugate to a non $\ZZ$ Seifert subgroup.

\noindent So that the last lemma allows one to decide for any
element $g$ in the group of such a triangulated 3-manifold whether
cases $(i)$, $(ii)$ or $(iii)$ occur and provide, if any, the
Seifert characteristic submanifold $S$ as well as a Seifert
subgroup $K=i_*(\pi_1(S))\subset\pi_1(M)$ containing an element
conjugate with $g$, together with a conjugating element.

\section{Proof of the main result}

We are now ready to give proofs of the results enonced in \S 1. We
are first concerned with the main result, that is theorem
\ref{main}. The following preliminary step reduces the proof to
the case of closed irreducible geometrizable 3-manifolds.

\begin{lem}\label{reduction}
Conjugacy problem in groups of non-oriented geometrizable
3-manifolds reduces to conjugacy problem in groups of closed
irreducible geometrizable 3-manifolds. Moreover if a triangulation
of $M$ as well as a solution to the word problem in $\pi_1(M)$ are
given, the reduction process applied to $\pi_1(M)$ is
constructive.
\end{lem}

\noindent{\sl Proof.} Let $M$ be a non-oriented geometrizable
3-manifold ; we are concerned with the conjugacy problem in
$\pi_1(M)$. We process the reduction in two steps.

Gluing a 3-ball to each spherical component of $\partial M$ leaves
$\pi_1(M)$ inchanged ; so we suppose in the following that $M$ has
no spherical boundary component. If $\partial M$ is non-empty,
double the manifold $M$ along its boundary to obtain the closed
non-oriented 3-manifold that we shall denote $2M$. Lemma 1.1 of
\cite{cp3mg} asserts that the inclusion map of $M$ in $2M$ induces
an embedding of $\pi_1(M)$ in $\pi_1(2M)$, and that
$u,v\in\pi_1(M)$ are conjugate in $\pi_1(M)$ if and only if they
are conjugate in $\pi_1(2M)$ ; hence the conjugacy problem in
$\pi_1(M)$ reduces to the conjugacy problem in $\pi_1(2M)$.
Moreover the closed manifold $2M$ is geometrizable ; indeed if
$\ol{M}$ and $\ol{2M}$ denote respectively the orientation covers
of $M$ and $2M$, one has that $\ol{2M}$ is the double of $\ol{M}$,
and since $M$ is geometrizable, lemma 1.2 of \cite{cp3mg} states
that $\ol{2M}$ is geometrizable. Hence the conjugacy problem in
$\pi_1(M)$ reduces to conjugacy problem in groups of closed
geometrizable 3-manifolds. Moreover  if $M$ is given by a
triangulation, the reduction can be achieved in a constructive way
for one can constructively produce a triangulation of $2M$ and the
natural embedding from $\pi_1(M)$ to $\pi_1(2M)$.

We are now concerned with the second step, and suppose $M$ to be
moreover closed. A Kneser-Milnor decomposition splits $M$ in a
connected sum of the prime closed geometrizable (non necessarily
non-oriented) factors $M_1, M_2,\ldots ,M_i$ and $\pi_1(M)$ splits
as the free product of the
$\pi_1(M_1),\pi_1(M_2),\ldots,\pi_1(M_n)$. Basic fact upon
conjugacy in free products (\cite{mks}) shows that the conjugacy
problem in $\pi_1(M)$ reduces to conjugacy problems in each of the
$\pi_1(M_i)$. Now either $M_i$ is a $S^2$-bundle over $S^1$, and
hence $\pi_1(M_i)=\Z$, or $M_i$ is irreducible. Hence conjugacy
problem in $\pi_1(M)$ reduces to conjugacy problem in groups of
closed irreducible geometrizable 3-manifolds. If $M$ is given by a
triangulation an algorithm in \cite{jt} for a Kneser-Milnor
decomposition allows to process the reduction in a constructive
way.\hfill$\square$\smallskip

We are now ready to show the main result upon conjugacy problem in
non oriented geometrizable manifolds.  Note that we not only prove
the existence of a solution : we rather show how, given a
triangulation of $M$, one can construct the algorithm.\\

\noindent {\bf Proof of theorem \ref{main}.} According to lemma
\ref{reduction}, and to a solution to the conjugacy problem in
groups of oriented geometrizable 3-manifolds (\cite{cp3mg}), we
are left with the case of groups of closed irreducible
non-oriented geometrizable 3-manifolds. In the following $M$
stands for a closed irreducible non-oriented geometrizable
3-manifold, and $p:N\longrightarrow M$ for the orientation cover
of $M$.
\begin{lem}\label{l4}
Given a triangulation of $M$ one can algorithmically produce a
triangulation of its orientation cover $N$ as well as the covering
map and covering automorphism.
\end{lem} \noindent{\sl Proof of
lemma \ref{l4}.} The triangulation of  $M$ can be easily given as
a triangulation of a $PL$-ball $B$ together with a gluing of pairs
of triangles in $\partial B$. Pick an orientation of $B$ ; it
induces an orientation of each triangle in $\partial B$. Identify
paired triangles in $\partial B$ each time their gluing preserves
orientation, to obtain a new oriented $PL$-manifold $C$ together
with orientation reversing gluing  of pairs of triangles in
$\partial C$. Consider a copy $C'$ of $C$ and give $C'$ the
reverse orientation. For each triangle $t$ in $\partial C$ denote
by $t'$ its copy in $\partial C'$. For each gluing of triangles
$t_1,t_2$ in $\partial C$, glue coherently in $C\cup C'$, $t_1$
with $t_2'$ and $t_1'$ with $t_2$. The manifold obtained is $N$
together with a triangulation, and the construction implicitly
produces the covering map $p:N\longrightarrow M$ as well a the
covering automorphism.\hfill $\square$\smallskip

We suppose in the following $M$ to be given by a triangulation.
The above lemma allows to produce a triangulation of $N$ as well
as the covering map $p:N\longrightarrow N$ and the covering
automorphism $\s$. Let $G=\pi_1(M)$ and $H=\pi_1(N)$.\smallskip

The oriented manifold $N$ may be reducible. In such a case (cf.
\cite{swarup}), $N$ contains a minimal system of essential
pairwise disjoint $\s$-invariant spheres $S=\{S_1,S_2,\ldots ,
S_q\}$, with
 image in $M$ a system of essential pairwise disjoint projective planes
$P=\{P_1,P_2,\ldots ,P_q\}$ such that : $(i)$ $N$ cutted along $S$
decomposes  into components $N_1,N_2,\ldots ,N_p$ and $n$
components homeomorphic to $S^2\times S^1$, such that each
manifold $\hat{N_i}$ obtained by gluing a ball to each $S^2\subset
\partial N_i$ is irreducible and non simply connected ; $(ii)$
$\pi_1(N)=\pi_1(N_1)*\cdots *\pi_1(N_p)*F(n)$ ; $(iii)$ $M$ cutted
along $P$ has component $M_1,M_2,\ldots, M_p$ as well as
components $\Bbb P_2\times I$, and the covering projection sends
$N_i$ onto $M_i$ and each $S^2\times I$ component onto a $\Bbb
P_2\times I$ component ; (iv) $\pi_1(M)$ splits as a graph of
group where vertex groups are $\pi_1(M_1),\ldots ,\pi_1(M_p)$ as
well as $n$ groups of order 2 and edge groups all have order 2.
Apply the following algorithm to find a system $S$ of pairwise
disjoint essential $\s$-invariant spheres.
\begin{lem}\label{l3}
One can algorithmically find a system of pairwise disjoint
essential $\s$-invariant spheres in $N$ and a system of pairwise
disjoint essential projective planes in $M$, as above.
\end{lem}
\noindent{\sl Proof of lemma \ref{l3}.} Apply in $N$ the algorithm
of \cite{jt} to find a system $S$ of disjoint spheres which
decomposes $N$ into pieces which, once balls are glued on the
boundary, are all irreducible. Apply the $S^3$ recognition
algorithm to detect all pieces in $N\setminus S$ homeomorphic to a
ball or  to $S^2\times I$. Delete each $S_1$ in $S$ which bounds
either a ball or two $S^2\times I$ pieces, each each couple
$S_1,S_2$ in $S$ whenever they are both separating and bound a
$S^2\times I$ piece. Up to this point $S$ consists only of
essential spheres. It remains to deform $S$ so that the spheres be
$\s$-invariant. For each $S_1\in S$, $\s(S_1)\cap
S_1\not=\emptyset$ cause otherwise $N$ wouldn't be irreducible.
For each $S_1\in S$ construct $\s(S_1)$, and use the $S^3$
recognition algorithm to construct an essential sphere included in
$S_1\cup\s(S_1)$ which is preserved under $\s$ ; change in $S$,
$S_1$ into this sphere. Each time two spheres in $S$ are not
disjoint apply an analog process to change them  into disjoint
essential spheres ; since the number of self intersection of
$S\cup\s(S)$ decreases the process must stop ; we are finally led
with a system of pairwise disjoint essential $\s$-invariant
spheres in $N$ ; their image under the covering projection gives
the system $P$ of essential pairwise disjoint projective planes in
$M$.\hfill $\square$\smallskip

Compute finite presentations of $G$ and $H$ in such a way that if
$S$ is non empty the generators of $H$ (respectively $G$) is the
union of the generators of the factors in its free product
decomposition (resp. generators of the vertex groups plus
additional stable letters). Compute the embedding of $H$ in $G$
(by the image of the generators of $H$ in $G$). Consider also
solutions to the word problem in $G$ and $H$, as well as a
solution to the conjugacy problem in $H$.
They can be given in a contructive way : an algorithm for  the
word problem in $G$ comes from the automatic group theory and can
be constructed (cf. \cite{epstein}) ; \cite{cp3mg} allows to
construct an algorithm for the conjugacy problem (and hence the
word problem) in $H$.

Use the solution to the word problem to decide whether $G$ is
abelian ; if so  a solution to the word problem gives a solution
to the conjugacy problem in $G$. So in the following we suppose
that $M$ is not homeomorphic to $\Bbb P_2\times S^1$.
\smallskip

 Suppose
$u$ and $v \in G$ are given and one wants to decide whether $u\sim
v$ in $G$. Find first respective classes of $u,v$ in $G/H$. This
can be done by using the naive algorithm of lemma \ref{ld3} or can
be achieved in a more efficient way by the following process :
suppose one knows for each generator $a$ of $G$ whether $a$ is
orientation reversing or not. Then $u$ lies in $H$ if and only if
the word  representing $u$ contains an even number of orientation
reversing generators. Note that from   a triangulation of $M$ one
can easily deduce a  generating set of $G$ and  check for each
generator whether it reverses orientation or not ; for example by
constructing a triangulation of $M$ as appearing in the proof of
lemma \ref{l4} : the gluing maps give generators of $G$.

 Since
$H$ has index 2 in $G$, if  $u$ and $v$ lie in distinct classes
they are definitely not conjugate in $G$. If $u$ and $v$ both lie
in $H$, then the lemma \ref{ld1} together with a solution to the
conjugacy problem in $H$ allow constructively to decide whether
$u$ and $v$ are conjugate in $G$. So in the following we will
suppose that $u$ and $v$ both lie in $G\setminus H$.

Use a solution to the word problem in $H$ to decide whether
$u^2=1$ and $v^2=1$. If exactly one of the relations occurs then
$u$ and $v$ are not conjugate in $G$, while if both relations
occur then the following lemma  allows to constructively decide
whether $u$ and $v$ are conjugate or not.

\begin{lem}
\label{ordre2} Let $G$ be as above. One can  decide for any couple
of order 2 elements $u, v\in G$ whether $u$ and $v$ are conjugate
in $G$.
\end{lem}
\noindent{\sl Proof.}  Consider the above system $P$ of essential
projective planes in $M$ and delete in $P$ one of $P_i,P_j$ every
time $P_i,P_j$ cobound a $\Bbb P_2\times I$ component. It follows
from  \cite{epsteinz2}, \cite{stalling}, \cite{swarup} that  each
order 2 element in $G$ is conjugate to some
$H_i=i_*(\pi_1(P_i))=\Z_2$ for $P_i\in P$, and moreover if we
denote $h_i$ the generator of $H_i$, the different $h_i$ have
disjoint conjugacy classes in $G$. Hence a naive algorithm goes as
follows : the different $h_i$ are given. Enumerate all the
elements $g$ of $G$ (as words on a given finite generating set),
and in parallel for each $g$ obtained and each $h_i$ use a
solution to the word problem in $G$ to decide whether $u=h_i^g$ or
$v=h_j^g$. The algorithm halts once it finds such $h_i,h_j$ in the
respective conjugacy classes of $u$ and $v$. Then $u\sim v$ if and
only if $h_i=h_j$. \hfill $\square$ \smallskip

 We will be concerned in the following with
the remaining case : $u,v$ both lie in $G\setminus H$ and both
have order different than 2 ; according to \cite{epsteinz2} both
$u$ and $v$ must have infinite order in $G$.

Use the lemma \ref{ld1} together with a solution to the conjugacy
problem in $H$ to decide whether $u^2\sim v^2$ in $G$. If it does
not arise then $u$ and $v$ are not conjugate in $G$. So we suppose
in the following that $u^2\sim v^2$ in $G$. Find an element $h\in
G$ such that $u^2=(v^2)^h$ in $G$. This can be naively performed
by enumerating all $g\in G$ and for each $g$ obtained by deciding
in parallel with a solution to the word problem whether
$u^2=(v^2)^g$ ; it can also be more efficiently achieved by going
into the line of the proofs of \cite{cp3mg} to remark that the
solution to the conjugacy problem
 provides such a conjugating element $h$.

 Up to this point the set $C_G(u^2,v^2)=h.Z_G(v^2)$ is non empty.
 It obviously
 contains the set
$C_G(u,v)$ ; note also that $Z_G(v^2)$ contains $Z_G(v)$ as a
subgroup as well as $Z_H(v^2)$ as an index 2 subgroup ; $Z_G(v^2)$
is generated by $Z_H(v^2)$ and $v$.\smallskip

We are now interested in the centralizer $Z_H(v^2)$. If the system
$S$ of essential spheres in $N$ is non-empty, $H$ splits non
trivially as a free product. If so make use of a solution to the
word problem in $H$ (together with elementary facts about free
products, cf. \cite{mks}) to write down respective cyclic
conjugates $dv^2d^{-1}$ of $v^2$ and $cu^2c^{-1}$ of $u^2$ in a
cyclically reduced form. If $dv^2d^{-1}$ has length greater than 1
(according to the splitting) then $Z_H(v^2)$ is infinite cyclic.
If $Z_H(v^2)$ is non cyclic $dv^2d^{-1}$ and $cu^2c^{-1}$ must
both lie in a factor $\pi_1(N_1)$ of the free product
decomposition of $H$.
We can apply in $N_1$  lemma \ref{central} together with its
following remark   to decide whether $dv^2d^{-1}$, and hence
$v^2$, has a cyclic or $\Z\oplus\Z$ centralizer, or other whether
$dv^2d^{-1}$ is conjugate to some non $\Z\oplus\Z$ Seifert
subgroup. Note that if $dv^2d^{-1}$ has a non cyclic centralizer,
it is conjugate to a Seifert subgroup $W$ of $\pi_1(N_1)$, and the
algorithm produces both $W$ and $e\in\pi_1(N_1)$ such that
$(ed)v^{2}(ed)^{-1}\in W$.
\smallskip

We consider first the case where $Z_H(v^2)$ is infinite cyclic ;
then $Z_G(v^2)$ contains $\Z$ as an index 2 subgroup. If
$Z_G(v^2)$ is torsion-free then it must be cyclic, say
$Z_G(v^2)=<w>$. But since $v\in Z_G(v^2)$, $v$ is a power of $w$,
and since $Z_G(v)\subset Z_G(v^2)$, it implies that
$Z_G(v)=Z_G(v^2)=<w>$. If $Z_G(v^2)$ has torsion, let us denote by
$t$ a generator of its index 2 subgroup $Z_H(v^2)$. The group
$Z_G(v^2)$ is  generated by $v$ and $t$ and must be one of the
groups appearing in the following lemma.
\begin{lem}\label{l1}
A torsion group $K$ with generators $v,t$, such that
$<t>\approx\Z$ has index 2 in $K$ must be one of :
$$<v,t\ |\ [v,t]=1, v^2=t^{2n}>\approx \Z\oplus\Z_2$$
$$<v,t\ |\ t^v=t^{-1}, v^2=1>\approx \Z_2 *\Z_2$$
\end{lem}
\noindent {\sl Proof of lemma \ref{l1}.}  The group $K$ admits the
presentation $<v,t\ |\ t^v=t^{\pm1}, v^2=t^p>$ with $p\in\Z$ ;
$K\setminus <t>$  contains an element $w$ with finite order
$k\not=0$. In particulary $w^k$ lies in $<t>$ so that $k$ must be
even ; hence $K$ contains an element $t^av$ with order 2.
Suppose first that $t^v=t$, so that $1=(t^{a}v)^2=t^{2a+p}$. It
follows that $p$ is even, say $p=2n$, which gives the first above
presentation.
Suppose then that $t^v=t^{-1}$ ; one has $1=(t^av)^2=v^2$ which
gives the second above presentation.\hfill $\square$\smallskip

The latter group cannot occur since $v$ has infinite order. For
the former group, since $[t,v]=1$, one has $Z_G(v)= Z_G(v^2)$.
Hence whenever  $Z_H(v^2)$ is infinite cyclic, $Z_G(v)=Z_G(v^2)$
and the lemma \ref{ld2} allows to decide whether $u\sim v$ in
$G$.\smallskip

We consider now the case where $v^2$ is conjugate to a Seifert
subgroup. We change $u,v,h$ respectively into $cuc^{-1}$,
$(ed)v(ed)^{-1}$ and $ch(ed)^{-1}$ ($c,d,e$ have been constructed
above), so that $u^2,v^2,h$ lie in the factor $\pi_1(N_1)$ of the
free product decomposition of $H$, and $u^2=(v^2)^h$. Moreover
$v^2$ lies in a given Seifert subgroup $W$ of $\pi_1(N_1)$.
 According to the following lemma we then have
$u,v\in\pi_1(M_1)$, where  $p(N_1)=M_1$.
\begin{lem}\label{l5} Suppose that $K$ splits as a fundamental group of a graph of group,
and that $xv^2x^{-1}$ lies in a vertex subgroup $G_V$ of $K$. Then
$xvx^{-1}$ also lies in $G_V$.
\end{lem}
\noindent{\sl Proof of lemma \ref{l5}.} Delete an edge in the
graph of group of $K$. It decomposes $K$ in an amalgam or an HNN
extension according to whether the edge is separating or not.
Denote by $K'$ a factor in this decomposition of $K$. If $w\in K$
is cyclically reduced and if $w^2$ lies in the factor $K'$ then
$w$ must also lie in $K'$ (cf. \cite{mks}, \cite{rotman}). If $w$
is non cyclically reduced, then it can be written in a reduced
form as $w=w_1\cdots w_k w' w_k^{-1}\cdots w_1^{-1}$ where $w'$ is
cyclically reduced. Hence $w^2=w_1\cdots w_k(w')^2w_k^{-1}\cdots
w_1^{-1}$ is reduced and non cyclically reduced and cannot lie in
the factor $K'$ of $K$. Note that $K'$ splits as a fundamental
group of a graph of group such that its vertex subgroups are
vertex subgroups of $K$. Pursue inductively the process by
decomposing $K'$ along an edge until we are led with factors which
are vertex subgroups of $K$. To conclude remark that if
$xv^2x^{-1}$ lies in a vertex subgroups $G_V$ of $K$ then, for any
of the above decompositions, its reduced form  is also cyclically
reduced and $xv^2x^{-1}$ must lie in some factor ;  the above
argument shows that the same conclusion applies to $xvx^{-1}$, so
that $xvx^{-1}$ must lie in $G_V$. \hfill$\square$\smallskip

Suppose first  that $v^2$ lies  in a $\Z\oplus\Z$ Seifert subgroup
of $\pi_1(N_1)$, so that $Z_H(v^2)=\Z\oplus\Z$.  The element $v$
obviously normalizes $Z_H(v^2)$ ; hence a $\Z\oplus\Z$ Seifert
submanifold homeomorphic to (an eventually punctured)
$T_1=S^1\times S^1\times I$ of $N_1$ can be chosen to be preserved
under the orientation reversing covering automorphism associated
with $v$ and its image in $M_1$ under the covering projection
induces the embedding of $Z_G(v^2)$ (generated by $Z_H(v^2)$ and
$v$) ; $Z_G(v^2)$ can only be the group of the Klein Bottle. One
must have :
$$Z_H(v^2)=<a,b\ |\ [a,b]=1>\quad Z_G(v^2)=<a,b,t\ |\ [a,b]=1, t^2=a, b^t=b^{-1}>$$
Moreover it's fairly easy to find generators $a,b,t$ of $Z_G(v^2)$
as above : pick an arbitrary base of
$Z_H(v^2)=i_*(\pi_1(T_1))\approx \Z\oplus\Z$, and using a solution
to the word problem  write down the conjugacy action of $v$ on
$Z_H(v^2)$ as a matrix $\cal{M}$ in $GL(2,\Z)$. Diagonalize
$\cal{M}$, it has eigen values $\l_1=1$ and $\l_2=-1$ ; $a,b$ are
the eigenvectors respectively associated with $\l_1$ and $\l_2$.
Then use a solution to the word problem in order to write down the
generator $t=a^nb^mv$.
\begin{lem}\label{l2}
Let $K=<a,b,t\ |\ [a,b]=1, t^2=a, b^t=b^{-1}>$ and
$v=a^{n_1}b^{n_2}t$. Then $Z_K(v)=<b^{n_2}t>\supset <a>$ and
$v'\sim v$ if and only if $v'=a^{n_1}b^{m_2}t$ with $m_2=n_2\mod
2$.
\end{lem}
\noindent {\sl Proof of lemma \ref{l2}.} Let $A$ be the abelian
subgroup generated by $a$ and $b$. Let $w$ be an element of $A$ ;
it's easy to see that if $w$ lies in $<a>$ then $Z_K(w)=K$ ($a$
and $t$ commute) and otherwise $Z_K(w)=A$. Let $v=a^{n_1}b^{n_2}t$
; with the above $Z_K(v)\cap A=<a>$. Let $z=a^{m_1}b^{m_2}t\in
K\setminus A$. Computation shows that
$vz=a^{n_1+m_1+1}b^{n_2-m_2}$ and $zv=a^{n_1+m_1+1}b^{m_2-n_2}$,
and hence $z\in Z_K(v)$ if and only if $m_2=n_2$ that is
$z=va^{m_1-n_1}$ for some $m_1\in \Z$. It follows that $Z_K(v)$ is
generated by $b^{n_2}t$ and $a$ ; but since $(b^{n_2}t)^2=a$,
$Z_K(v)=<b^{n_2}t>$.
We are now concerned with the conjugacy class of $v$. One has that
$a$ and $v$ commute, and for $n\in\Z$,
$b^nvb^{-n}=a^{n_1}b^{n_2+2n}t$, while $t^nvt^{-n}$ equals $v$
when $n$ is even and $a^{n_1}b^{-n_2}t$ otherwise. Hence $v'$ and
$v$ are conjugate if and only if $v'=a^{n_1}b^{m_2}t$ with
$m_2=n_2\mod 2$.\hfill$\square$\smallskip

 Write down $v$ on the generators $a,b,t$, say $v=a^{n_1}b^{n_2}t$
 ; with the above lemma
$Z_G(v)=<b^{n_2}t>\supset <a>$, and $v'$ is conjugate in
$Z_G(v^2)$ to $v$ if and only if $v'=a^{n_1}b^{m_2}t$ with
$m_2=n_2 \mod 2$.

The set $C_G(u^2,v^2)=hZ_G(v^2)$ is given and contains $C_G(u,v)$.
Hence $u$ and $v$ are conjugate in $G$ if and only if there exists
an element $w\in Z_G(v^2)\setminus Z_G(v)$ such that
$h^{-1}uh=wvw^{-1}$. With the above, $u\sim v$ if and only if
$h^{-1}uh=a^{n_1}b^{m_2}t$ with $m_2=n_2\mod 2$, that is if and
only if $h^{-1}uh$ lies in $v<b^2>$. Use  a solution to the word
problem in $G$ to decide first whether $h^{-1}uh$ lies in
$Z_G(v^2)$, and if yes to secondly write $h^{-1}uhv^{-1}$ on the
generators $a,b,t$ and then check whether $h^{-1}uhv^{-1}$ lies in
$<b^2>$. This process allows to decide whether $u$ and $v$ are
conjugate in $G$ other not.\smallskip

Finally suppose that $v^2$ lies in a non $\ZZ$ Seifert subgroup
$W$ of $H$. One has splittings of $\pi_1(N_1)$ and $\pi_1(M_1)$ as
fundamental groups of graphs of group induced by the
JSJ-decompositions of the irreducible manifold $\hat{N_1}$ and
$M_1$, so that vertex groups of $\pi_1(N_1)$ are index 1 or 2
subgroups of vertex groups of $\pi_1(M_1)$.
 The element $v^2$ lies in a vertex subgroup $H_s=W$ coming
from a Seifert piece $N'$ in the JSJ splitting of $\hat{N_1}$.
Change $u$ into $huh^{-1}$ so that $u^2=v^2$. According to the
lemma \ref{l5}, $u,v$ lie in a vertex subgroup $G_s=p_*(H_s)$ of
$\pi_1(M_1)$, such that $G_s\cap H=H_s$ is an index 2 subgroup of
$G_s$ (since $v\not\in H$). In fact $G_s=\pi_1(M')$ where $M'$ is
a piece in the JSJ decomposition of $M_1$ and one has the
orientation cover $p':N'\setminus \bigcup B^2\longrightarrow M'$ ;
the covering automorphism can be extended to an orientation
reversing involution of $N'$ with a finite number of fixed points.

  On the one hand  $C_H(u^2,v^2)=Z_H(v^2)$ is included in $H_s$ (cf. theorem VI.I.6,
  \cite{js} or remark \S 3), and on the other $C_G(u^2,v^2)=Z_G(v^2)$ is
  generated by $Z_H(v^2)$ and $v$, and hence is included in $G_s$. Since
  $C_G(u^2,v^2)\supset C_G(u,v)$, if $u$ and $v$ are conjugate in
  $G$, they must be conjugate in $G_s$.

  The manifold $N'$ is a Seifert fibered space which cannot be
  modelled on $Nil$ geometry (such manifolds do not admit an
  orientation reversing involution, cf. \cite{scott}) and
  according to \cite{nr1} $H_s$ is a biautomatic group. Hence (cf.
  \cite{epstein}) $G_s$ is also biautomatic and one can construct
  an algorithm to solve the conjugacy problem in $G_s$. So that we
  can decide whether $u$ and $v$ are conjugate in $G$ or not.\hfill
  $\square$\medskip

\noindent {\bf Proof of theorem \ref{t2}.} Double the 3-manifold
$M$ along the identity on $F$ ; one obtains a 3-manifold $2M$. The
argument in the proof of lemma 1.2 of \cite{cp3mg} as well as the
observation that the orientation cover of $2M$ is the double of
the orientation cover of $M$ along lifting(s) of $F$ show that
$2M$ is geometric. Its group splits into an amalgam of two copies
of $M$ along $H$, say $\G=G*_H G$ ; given $g$ in the $G$ left
factor we note $\bar{g}$ the corresponding element of the $G$
right factor. Since the gluing is along the identity $h=\bar{h}$,
one has that $u,\bar{u}\in G$ are equal (resp. conjugate) in $G$
if and only if $u\in H$ (resp. $u$ conjugate in $G$ to some $h\in
H$). Hence a solution to the word problem (resp. conjugacy
problem) in $\pi_1(2M)$ provides the algorithm.\hfill $\square$

\end{document}